\documentclass[reqno]{amsart}
\usepackage{amssymb,mathtools,anysize}
\marginsize{2.5cm}{2.5cm}{2cm}{2cm}
\usepackage{ifpdf}
\ifpdf
\usepackage[hyperindex]{hyperref}%,pagebackref
\else
\expandafter\ifx\csname dvipdfm\endcsname\relax
\usepackage[hypertex,hyperindex]{hyperref}
\else
\usepackage[dvipdfm,hyperindex]{hyperref}
\fi
\fi
\allowdisplaybreaks[4]
\numberwithin{equation}{section}
\theoremstyle{plain}
\newtheorem{thm}{Theorem}[section]

\theoremstyle{definition}
\newtheorem{dfn}{Definition}[section]
\theoremstyle{remark}
\newtheorem{rem}{Remark}[section]
\DeclareMathOperator{\td}{\textup{d}}

\newcommand{\cmdeg}[1]{\sideset{}{_\textup{cm}^{#1}}\deg}

\begin{document}

\title[Completely monotonic degrees of digamma function]
{Partial solutions to several conjectures on completely monotonic degrees for remainders in asymptotic expansions of the digamma function}

\author[F. Qi]{Feng Qi}
\address[Qi]{Institute of Mathematics, Henan Polytechnic University, Jiaozuo 454010, Henan, China; School of Mathematical Sciences, Tianjin Polytechnic University, Tianjin 300387, China}
\email{\href{mailto: F. Qi <qifeng618@gmail.com>}{qifeng618@gmail.com}, \href{mailto: F. Qi <qifeng618@hotmail.com>}{qifeng618@hotmail.com}, \href{mailto: F. Qi <qifeng618@qq.com>}{qifeng618@qq.com}}
\urladdr{\url{https://qifeng618.wordpress.com}, \url{https://orcid.org/0000-0001-6239-2968}}

\author[M. Mahmoud]{Mansour Mahmoud}
\address[Mahmoud]{Mathematics Department, Faculty of Science, Mansoura University, Mansoura 35516, Egypt}
\email{mansour@mans.edu.eg}
\urladdr{\url{https://orcid.org/0000-0002-5918-1913}}

\dedicatory{Dedicated to people facing and battling COVID-19}

\begin{abstract}
Motivated by several conjectures posed in the paper ``F. Qi and A.-Q. Liu, \textit{Completely monotonic degrees for a difference between the logarithmic and psi functions}, J. Comput. Appl. Math. \textbf{361} (2019), 366\nobreakdash--371; available online at \url{https://doi.org/10.1016/j.cam.2019.05.001}'', the authors compute several completely monotonic degrees of the remainders in the asymptotic expansions of the logarithm of the gamma function and in the asymptotic expansions of the logarithm of the digamma function.
\end{abstract}



\thanks{This paper was typeset using\AmS-\LaTeX}

\maketitle
\tableofcontents

\section{Preliminaries}

The classical Euler gamma function $\Gamma(z)$ can be defined by
\begin{equation*}
\Gamma(z)=\int_0^{\infty}t^{z-1}e^{-t}\td t, \quad\Re(z)>0.
\end{equation*}
The logarithmic derivative of $\Gamma(x)$, denoted by $\psi(x)=\frac{\Gamma'(x)}{\Gamma(x)}$, is called the digamma function and the derivatives $\psi^{(i)}(x)$ for $i\ge0$ are called the polygamma functions. For more information on $\Gamma(z)$ and $\psi^{(i)}(x)$, refer to~\cite{abram, NIST-HB-2010, singularity-combined.tex, Yang-Tian-JMAA-2019} and closely related references therein.
\par
Recall from~\cite[Chapter~1]{Schilling-Song-Vondracek-2nd} that, if a function $h$ on an interval $I$ has derivatives of all orders on $I$ and
\begin{equation*}
(-1)^nh^{(n)}(t)\geq0, \quad t\in I, \quad n\in\{0\}\cup\mathbb{N},
\end{equation*}
then we call $h$ a completely monotonic function on $I$. In other words, a function $h$ is completely monotonic on an interval $I$ if its odd derivatives are negative and its even derivatives are positive on $I$.
Theorem~12b in~\cite[p.~161]{widder} states that a necessary and sufficient condition for $h$ to be completely monotonic on the infinite interval $(0,\infty)$ is that
\begin{equation}\label{Theorem12b-Laplace}
h(t)=\int_0^\infty e^{-ts}\td\sigma(s), \quad s\in(0,\infty),
\end{equation}
where $\sigma(s)$ is non-decreasing and the above integral converges for $s\in(0,\infty)$. In other words, a function is completely monotonic on $(0,\infty)$ if and only if it is a Laplace transform of a non-negative measure.

\begin{dfn}\label{cm-degree-def}
Let $h(t)$ be a completely monotonic function on $(0,\infty)$ and let $h(\infty)=\lim_{t\to\infty}h(t)\ge0$.
If the function $t^\alpha[h(t)-h(\infty)]$ is completely monotonic on $(0,\infty)$ if and only if $0\le\alpha\le r\in\mathbb{R}$, then we say that $h(t)$ is of completely monotonic degree $r$; if $t^\alpha[h(t)-h(\infty)]$ is completely monotonic on $(0,\infty)$ for all $\alpha\in\mathbb{R}$, then we say that the completely monotonic degree of $h(t)$ is $\infty$.
\end{dfn}

The function $t^\alpha[h(t)-h(\infty)]$ in Definition~\ref{cm-degree-def} can be essentially regarded as the ratio $\frac{h(t)-h(\infty)}{1/t^\alpha}$ between the completely monotonic function $h(t)-h(\infty)$ and the $\alpha$ power of the completely monotonic function $\frac{1}{t}$ on $(0,\infty)$. This is the reason why we designed in~\cite{psi-proper-fraction-degree-two.tex} a notation $\cmdeg{t}[h(t)]$ to denote the completely monotonic degree $r$ of $h(t)$ with respect to $t\in(0,\infty)$. According to this idea, we can define the completely monotonic degree $r$ of $h(t)$ with respect to $g(t)$, denoted by $\cmdeg{g(t)}[h(t)]$, as the largest number $\alpha$ such that the ratio $\frac{h(t)-h(\infty)}{[g(t)-g(\infty)]^\alpha}$ between the completely monotonic function $h(t)$ and the $\alpha$ power of the completely monotonic function $g(t)$ on $(0,\infty)$.
\par
Proposition~1.2 in~\cite{22ICFIDCAA-Filomat.tex} can be modified as that $\cmdeg{t}[h(t)]=r>0$ if and only if
\begin{equation*}
h(t)-h(\infty)=\int_0^\infty\biggl[\frac1{\Gamma(r)} \int_0^s(s-\tau)^{r-1}\td\mu(\tau)\biggr]e^{-ts}\td s
\end{equation*}
for $0<t<\infty$, where $\mu(\tau)$ is a bounded and non-decreasing measure on $(0,\infty)$.
For more information on completely monotonic degrees and their properties, please refer to the papers~\cite{Koumandos-MC-08, Koumandos-Lamprecht-MC-2010, Koumandos-Lamprecht-MC-2013, Koumandos-Pedersen-Monatsh-2011, Koumandos-Pedersen-09-JMAA, AIMS-Math-2019595.tex, Qi-Agar-Surv-JIA.tex} and closely related references therein.

\section{Motivations and main results}

In~\cite[pp.~374--375, Theorem~1]{psi-alzer} and~\cite[Theorem~1]{theta-new-proof.tex-BKMS}, the function
\begin{equation*}
t^\alpha[\ln t-\psi(t)], \quad\alpha\in\mathbb{R}
\end{equation*}
was proved to be completely monotonic on $(0,\infty)$ if and only if $\alpha\le1$. This means that completely monotonic degree of $\ln t-\psi(t)$ on $(0,\infty)$ is
\begin{equation}\label{psi-ln-degree=1}
\cmdeg{t}[\ln t-\psi(t)]=1.
\end{equation}
\par
In~\cite[Theorem~1.7]{10.1.1.97.1017.pdf}, the function
\begin{equation*}
t^2[\psi(t)-\ln t]+\frac{t}{2}
\end{equation*}
was proved to be decreasing and convex on $(0,\infty)$ and, as $t\to\infty$, to tend to $-\frac{1}{12}$.
In~\cite[Theorem~1]{Chen-Qi-Srivastava-09.tex}, the function
\begin{equation*}
t^2[\psi(t)-\ln t]+\frac{t}{2}+\frac{1}{12}
\end{equation*}
was verified to be completely monotonic on $(0,\infty)$. In~\cite[Theorem~2]{CAM-D-18-02975.tex}, the completely monotonic degree of
\begin{equation*}
\phi(t)=\psi(t)-\ln t+\frac{1}{2t}+\frac{1}{12t^2}
\end{equation*}
with respect to $t\in(0,\infty)$ was proved to be
\begin{equation}\label{phi-degree=2}
\cmdeg{t}[\phi(t)]=2.
\end{equation}

In~\cite[Theorem~8]{psi-alzer}, \cite[Theorem~2]{Koumandos-jmaa-06}, and~\cite{Xu-Han-SM-2009}, the functions
\begin{equation*}
R_n(t)=(-1)^n\Biggl[\ln\Gamma(t)-\biggl(t-\frac12\biggr)\ln t+t-\frac12\ln(2\pi)-\sum_{k=1}^n\frac{B_{2k}}{(2k)(2k-1)}\frac1{t^{2k-1}}\Biggr]
\end{equation*}
for $n\ge0$ were proved to be completely monotonic on $(0,\infty)$, where an empty sum is understood to be $0$ and the Bernoulli numbers $B_n$ are generated~\cite{CAM-D-18-00067.tex, 2Closed-Bern-Polyn2.tex} by
\begin{equation}\label{Bernoullu=No-dfn-Eq}
\frac{z}{e^z-1}=\sum_{n=0}^\infty B_n\frac{z^n}{n!}=1-\frac{z}2+\sum_{k=1}^\infty B_{2k}\frac{z^{2k}}{(2k)!}, \quad\vert z\vert<2\pi.
\end{equation}
This conclusion implies that the functions $(-1)^mR_n^{(m)}(t)$ for $m,n\ge0$ are completely monotonic on $(0,\infty)$. See also~\cite[Section~1.4]{Sharp-Ineq-Polygamma-Slovaca.tex} and~\cite[Theorem~3.1]{Mortici-AA-10-134}. By the way, we call the function $(-1)^nR_n(t)$ for $n\ge0$ remainders of the asymptotic formula for $\ln\Gamma(t)$. See~\cite[p.~257, 6.1.40]{abram} and~\cite[p.~140, 5.11.1]{NIST-HB-2010}.
\par
In~\cite[Theorem~1]{psi-alzer}, \cite[Theorem~1]{theta-new-proof.tex-BKMS}, and~\cite[Theorem~3]{Xu=Cen-Qi-Conj-JIA2020}, the degree $\cmdeg{t}[-R_0'(t)]=1$ was verified once again.
\par
In~\cite[Theorem~2.1]{Koumandos-Pedersen-09-JMAA}, it was proved that
\begin{equation}\label{Koumandos-Pedersen-09-JMAA-thm2.1}
\cmdeg{t}\bigl[R_n(t)\bigr]\ge n, \quad n\ge0.
\end{equation}
\par
In~\cite[Theorem~1]{Chen-Qi-Srivastava-09.tex}, \cite[Theorem~2]{CAM-D-18-02975.tex}, and~\cite[Theorem~4]{Xu=Cen-Qi-Conj-JIA2020}, the degree $\cmdeg{t}[-R_1'(t)]=2$ was proved once again.
\par
In~\cite[Theorems~1 and~2]{Xu=Cen-Qi-Conj-JIA2020}, it was shown that
\begin{equation*}
\cmdeg{t}\bigl[(-1)^2R_0''(t)\bigr]=2 \quad\text{and}\quad
\cmdeg{t}\bigl[(-1)^2R_1''(t)\bigr]=3.
\end{equation*}
\par
In~\cite{CM-deg-tri-p5.tex, CM-degree-tri-p5.tex}, Qi proved that
\begin{equation*}
4\le\cmdeg{t}\bigl[(-1)^2R_2''(t)\bigr]\le5.
\end{equation*}
\par
Due to the above results, we modify Qi's conjectures posed in~\cite{CAM-D-18-02975.tex} as follows:
\begin{enumerate}
\item
the completely monotonic degrees of $R_n(t)$ for $n\ge0$ with respect to $t\in(0,\infty)$ satisfy
\begin{equation}\label{R-01-0th-deg}
\cmdeg{t}[R_0(t)]=0, \quad\cmdeg{t}[R_1(t)]=1,
\end{equation}
and
\begin{equation}\label{R-n-0th-deg}
\cmdeg{t}[R_n(t)]=2(n-1), \quad n\ge2;
\end{equation}
\item
the completely monotonic degrees of $-R_n'(t)$ for $n\ge0$ with respect to $t\in(0,\infty)$ satisfy
\begin{equation}\label{Qi-Conj-Gam-Deg-der1-1-2}
\cmdeg{t}[-R_0'(t)]=1, \quad\cmdeg{t}[-R_1'(t)]=2,
\end{equation}
and
\begin{equation}\label{Qi-Conj-Gam-Deg-der1(2n-1)}
\cmdeg{t}[-R_n'(t)]=2n-1, \quad n\ge2;
\end{equation}
\item
the completely monotonic degrees of $(-1)^mR_n^{(m)}(t)$ for $m\ge2$ and $n\ge0$ with respect to $t\in(0,\infty)$ satisfy
\begin{equation}\label{R-01-m-th-deg}
\cmdeg{t}\bigl[(-1)^mR_0^{(m)}(t)\bigr]=m, \quad\cmdeg{t}\bigl[(-1)^mR_1^{(m)}(t)\bigr]=m+1,
\end{equation}
and
\begin{equation}\label{R-n-m-th-deg}
\cmdeg{t}\bigl[(-1)^mR_n^{(m)}(t)\bigr]=m+2(n-1), \quad n\ge2.
\end{equation}
\end{enumerate}
\par
In this paper, we will confirm that Qi's conjectures expressed in~\eqref{R-01-0th-deg} and~\eqref{Qi-Conj-Gam-Deg-der1-1-2} are true and, via a double inequality, partially confirm that the conjecture expressed in~\eqref{Qi-Conj-Gam-Deg-der1(2n-1)} is almost true. Our main results can be stated as the following theorem.

\begin{thm}\label{1der-deg=2n-1thm}
For $n\ge0$, the completely monotonic degrees of the remainder $R_n(t)$ with respect to $t\in(0,\infty)$ satisfy the two equalities in~\eqref{R-01-0th-deg}.
\par
For $n\ge0$, the completely monotonic degrees of the functions
\begin{equation*}
-R_n'(t)=(-1)^{n+1}\Biggl[\psi(t)-\ln t+\frac{1}{2t}+\sum_{k=1}^n\frac{B_{2k}}{2k}\frac1{t^{2k}}\Biggr]
\end{equation*}
with respect to $t\in(0,\infty)$ satisfy the two equalities in~\eqref{Qi-Conj-Gam-Deg-der1-1-2} and a double inequality
\begin{equation}\label{Qi-Conj-Gam-Deg-der(2n-1)-2n}
2n-1\le\cmdeg{t}[-R_n'(t)]\le2n, \quad n\ge2.
\end{equation}
\end{thm}

\section{Proof of Theorem~\ref{1der-deg=2n-1thm}}

Now we start out to prove our main results stated in Theorem~\ref{1der-deg=2n-1thm}.

\subsection{Proofs of equalities in~\eqref{R-01-0th-deg}}%
From Binet's first formula
\begin{equation*}
\ln\Gamma(t)= \biggl(t-\frac{1}{2}\biggr)\ln t-t+\ln \sqrt{2\pi}\,+\int_{0}^{\infty}\biggl(\frac{1}{e^{u}-1} -\frac{1}{u}+\frac{1}{2}\biggr)\frac{e^{-tu}}{u}\td u
\end{equation*}
for $t>0$ in~\cite[p.~28, Theorem~1.6.3]{andrews-askey-roy-1999b} and~\cite[p.~11]{magnus}, it is easy to see that
\begin{equation*}
R_0(t)=\ln\Gamma(t)-\biggl(t-\frac12\biggr)\ln t+t-\frac12\ln(2\pi)
\end{equation*}
and
\begin{equation*}
R_1(t)=-\biggl[\ln\Gamma(t)-\biggl(t-\frac12\biggr)\ln t+t-\frac12\ln(2\pi)-\frac1{12t}\biggr]
\end{equation*}
satisfy
\begin{equation*}
\lim_{t\to\infty}R_0(t)=\lim_{t\to\infty}R_1(t)=0.
\end{equation*}
The inequality~\eqref{Koumandos-Pedersen-09-JMAA-thm2.1} means that the completely monotonic degree of the remainder $R_n(t)$ for $n\ge0$ with respect to $t\in(0,\infty)$ is at least $n$. In particular, we have
\begin{equation}\label{R01-deg-ge01}
\deg^{t}[R_0(t)]\ge0 \quad\text{and}\quad \deg^{t}[R_1(t)]\ge1.
\end{equation}
If $t^\theta R_0(t)$ and $t^\lambda R_1(t)$ were completely monotonic on $(0,\infty)$, then their first derivatives should be non-positive. As a result, using
\begin{equation*}
R_0'(t)=\psi(t)-\ln t+\frac{1}{2t}\quad\text{and}\quad
-R_1'(t)=\psi(t)-\ln t+\frac{1}{2t}+\frac{1}{12t^2},
\end{equation*}
we acquire
\begin{equation*}
\theta\le-\frac{tR_0'(t)}{R_0(t)}
=\frac{t\bigl[\ln t-\psi(t)]-\frac{1}{2}}{\ln\Gamma(t)-\bigl(t-\frac12\bigr)\ln t+t-\frac12\ln(2\pi)}
\to0, \quad t\to0^+
\end{equation*}
and
\begin{equation*}
\lambda\le-\frac{tR_1'(t)}{R_1(t)}
=-\frac{t^2\bigl[\psi(t)-\ln t+\frac{1}{2t}+\frac1{12t^{2}}\bigr]} {t\bigl[\ln\Gamma(t)-\bigl(t-\frac12\bigr)\ln t+t-\frac12\ln(2\pi)-\frac1{12t}\bigr]}\to1, \quad t\to0^+,
\end{equation*}
where we used the limits
\begin{gather*}
\begin{aligned}
\lim_{t\to0^+}(t[\ln t-\psi(t)])&=1,\quad \text{(see~\cite[p.~374]{psi-alzer})}\\
\lim_{t\to0^+}\biggl(t^2\biggl[\psi(t)-\ln t+\frac{1}{2t}+\frac1{12t^{2}}\biggr]\biggr)&=\frac{1}{12}, \quad \text{(see~\cite[Theorem~1]{CAM-D-18-02975.tex})}
\end{aligned}\\
\begin{aligned}
R_0(t)&=\ln\frac{\Gamma(t+3)}{(t+2)(t+1)t}-\biggl(t-\frac12\biggr)\ln t+t-\frac12\ln(2\pi)\\
&=\ln\Gamma(t+3)-\ln(t+2)-\ln(t+1)-t\ln t-\frac12\ln t+t-\frac12\ln(2\pi)\\
&\to\infty,\quad t\to0^+,
\end{aligned}
\end{gather*}
and
\begin{align*}
tR_1(t)&=-t\biggl[R_0(t)-\frac1{12t}\biggr]\\
&=\frac{1}{12}-t\ln\Gamma(t+3)+t\ln\frac{t+2}{t+1}+\biggl(t+\frac{1}{2}\biggr)t\ln t-t^2+\frac{\ln(2\pi)}{2}t\\
&\to\frac{1}{12}, \quad t\to0^+.
\end{align*}
Consequently, it follows that
\begin{equation}\label{R01-deg-le01}
\deg^{t}[R_0(t)]\le0 \quad\text{and}\quad \deg^{t}[R_1(t)]\le1.
\end{equation}
Combining those inequalities in~\eqref{R01-deg-ge01} and~\eqref{R01-deg-le01} concludes those equalities in~\eqref{R-01-0th-deg}.

\subsection{Proofs of equalities in~\eqref{Qi-Conj-Gam-Deg-der1-1-2}}%
Since $-R_1'(t)=\phi(t)$,
\begin{gather*}
\lim_{t\to\infty}t[\ln t-\psi(t)]=\frac12,\quad \text{(see~\cite[p.~374]{psi-alzer})}\\
\lim_{t\to0^+}\frac{t\bigl[\frac{1}{2t^2}+\frac{1}{t}-\psi'(t)\bigr]}{\frac{1}{2 t}-\ln t+\psi(t)}
=\lim_{t\to0^+}\frac{\frac{1}{2}+t-t^2\bigl[\psi'(t+1)+\frac{1}{t^2}\bigr]}{\frac{1}{2}-t[\ln t-\psi(t)]}=1,
\end{gather*}
by the equation~\eqref{psi-ln-degree=1}, we conclude the first conclusion in~\eqref{Qi-Conj-Gam-Deg-der1-1-2}.
\par
The equation~\eqref{phi-degree=2} established in~\cite[Theorem~2]{CAM-D-18-02975.tex} is equivalent to the second conclusion in~\eqref{Qi-Conj-Gam-Deg-der1-1-2}.

\subsection{Proof of the equality~\eqref{Qi-Conj-Gam-Deg-der1(2n-1)}}%
Let
\begin{equation*}
f_n(v)= (-1)^n\Biggl[\frac{1}{v}-\frac{1}{2}\coth\frac{v}{2} +\sum_{k=1}^n\frac{B_{2k}}{(2k)!} v^{2k-1}\Biggr]
\end{equation*}
for $n\ge0$, where the empty sum is understood to be $0$.
Then, by virtue of the formulas
\begin{gather*}
\coth v=\frac{1+e^{-2v}}{1-e^{-2v}}=\frac{2}{1-e^{-2v}}-1,\\
\psi(z)=\ln z+\int_{0}^{\infty}\biggl(\frac{1}{v}-\frac{1}{1-e^{-v}}\biggr)e^{-zv}\td v, \quad \Re(z)>0
\end{gather*}
in~\cite[p.~140, 5.9.13]{NIST-HB-2010}, and
\begin{equation*}
\frac1{z^w}=\frac1{\Gamma(w)}\int_0^\infty v^{w-1}e^{-zv}\td v, \quad \Re(z), \Re(w)>0
\end{equation*}
in~\cite[p.~255, 6.1.1]{abram}, we derive that
\begin{align*}
(-1)^n\int_0^{\infty} f_n(v) e^{-tv}\td v
&=\int_0^{\infty}\Biggl[\frac{1}{v}-\frac{1}{2}\coth\frac{v}{2} +\sum_{k=1}^n\frac{B_{2k}}{(2k)!} v^{2k-1}\Biggr]e^{-tv}\td v\\
&=\int_0^{\infty}\Biggl[\frac{1}{v}-\frac1{1-e^{-v}}+\frac12 +\sum_{k=1}^n\frac{B_{2k}}{(2k)!} v^{2k-1}\Biggr]e^{-tv}\td v\\
&=\psi(t)-\ln t+\frac1{2t}+\sum_{k=1}^n\frac{B_{2k}}{2k}\frac1{t^{2k}}
\end{align*}
for $n\ge0$.
This means that
\begin{equation}\label{Remainder-f-n-Rel-Eq}
R_n'(t)=\int_0^{\infty} f_n(v) e^{-tv}\td v, \quad n\ge0.
\end{equation}
\par
The inequality~\eqref{Koumandos-Pedersen-09-JMAA-thm2.1} tells us that the remainder $R_n(t)$ for $n\ge0$ is completely monotonic in $t\in(0,\infty)$.
Then, by definition of completely monotonic functions, it is ready that the remainder $-R_n'(t)$ for $n\ge0$ is completely monotonic on $(0,\infty)$.
\par
Since the function $\frac{v}{e^v-1}-1+\frac{v}2$ is even in $v\in\mathbb{R}$, by virtue of the equation~\eqref{Bernoullu=No-dfn-Eq}, it follows that
\begin{align*}
f_n(v)
&=\frac{(-1)^n}v\Biggl[1-\frac{v}{1-e^{-v}}+\frac{v}2 +\sum_{k=1}^n\frac{B_{2k}}{(2k)!} v^{2k}\Biggr]\\
&=\frac{(-1)^n}v\Biggl[-\biggl(\frac{v}{e^v-1}-1+\frac{v}2\biggr) +\sum_{k=1}^n\frac{B_{2k}}{(2k)!} v^{2k}\Biggr]\\
&=(-1)^{n+1}\sum_{k=n+1}^{\infty}\frac{B_{2k}}{(2k)!} v^{2k-1}
\end{align*}
for $|v|<2\pi$ and $n\ge0$.
Accordingly, we have
\begin{equation}\label{f(n)2zero=zero}
\lim_{v\to0}f_n^{(\ell)}(v)=(-1)^{n+1}\lim_{v\to0} \sum_{k=n+1}^{\infty}\frac{B_{2k}}{(2k)!} \langle2k-1\rangle_\ell v^{2k-\ell-1}=0
\end{equation}
for $0\le\ell\le2n$ or $\ell=2m$ with $m,n\ge0$, where
\begin{equation*}
\langle \alpha\rangle_n=
\prod_{k=0}^{n-1}(\alpha-k)=
\begin{cases}
\alpha(\alpha-1)\dotsm(\alpha-n+1), & n\ge1\\
1,& n=0
\end{cases}
\end{equation*}
is called the falling factorial of $\alpha$. For more information on the falling and rising factorials, refer to the papers~\cite{CDM-68111.tex}.
\par
Recall from~\cite[Theorem~2.1]{Eight-Identy-More.tex}, \cite[Theorem~2.1]{exp-derivative-sum-Combined.tex}, and~\cite[Theorem~3.1]{CAM-D-13-01430-Xu-Cen} that, when $\vartheta>0$ and $t\ne-\frac{\ln\vartheta}\theta$ or when $\vartheta<0$ and $t\in\mathbb{R}$, we have
\begin{equation}\label{id-gen-new-form1}
\frac{\td^k}{\td t^k}\biggl(\frac1{\vartheta e^{\theta t}-1}\biggr)
=(-1)^k\theta^k\sum_{p=1}^{k+1}{(p-1)!S(k+1,p)}\biggl(\frac1{\vartheta e^{\theta t}-1}\biggr)^p
\end{equation}
for $k\ge0$, where
\begin{equation*}
S(k,p)=\frac1{p!}\sum_{q=1}^p(-1)^{p-q}\binom{p}{q}q^{k}, \quad 1\le p\le k
\end{equation*}
are the Stirling numbers of the second kind. Taking $\vartheta=\theta=1$ in~\eqref{id-gen-new-form1} leads to
\begin{equation}\label{id-gen-new-theta=1}
\frac{\td^k}{\td t^k}\biggl(\frac1{e^t-1}\biggr)
=(-1)^k\sum_{p=1}^{k+1}{(p-1)!S(k+1,p)}\biggl(\frac1{e^t-1}\biggr)^p
\end{equation}
for $k\ge0$. Utilizing~\eqref{id-gen-new-theta=1} results in
\begin{equation}\label{coth-infinity=0}
\lim_{v\to\infty}\frac{\td^\ell}{\td v^\ell}\biggl(\frac{1}{v}-\frac{1}{2}\coth\frac{v}{2}\biggr)
=\lim_{v\to\infty}\frac{\td^\ell}{\td v^\ell}\biggl[\frac{1}{v}-\biggl(\frac{1}{e^v-1}+\frac{1}{2}\biggr)\biggr]
=\begin{cases}0, & \ell\ge1;\\
-\dfrac12, & \ell=0.
\end{cases}
\end{equation}
Making use of~\eqref{coth-infinity=0} and
\begin{equation*}
\lim_{v\to\infty}\bigl(v^m e^{-tv}\bigr)=0, \quad t>0, \quad m\ge0
\end{equation*}
yields
\begin{equation}\label{f(n)2inifnity=zero}
\lim_{v\to\infty}\bigl[f_n^{(\ell)}(v)e^{-tv}\bigr]
=(-1)^n\lim_{v\to\infty}\Biggl[\biggl(\frac{1}{v}-\frac{1}{2}\coth\frac{v}{2}\biggr)^{(\ell)}e^{-tv}
+\sum_{k=1}^n\frac{B_{2k}}{(2k)!} \langle2k-1\rangle_\ell v^{2k-\ell-1}e^{-tv}\Biggr]
=0
\end{equation}
for $\ell,n\ge0$ and $t>0$.
Consequently, by the limits~\eqref{f(n)2zero=zero} and~\eqref{f(n)2inifnity=zero}, integrating by parts inductively $2n-1$ times in the equation~\eqref{Remainder-f-n-Rel-Eq} gives
\begin{align*}
R_n'(t)&=-\frac{1}{t}\int_0^{\infty} f_n(v) \td e^{-tv}\\
&=-\frac1t\biggl[f_n(v)e^{-tv}\bigl|_{v=0}^\infty-\int_0^{\infty} f_n'(v) e^{-tv}\td v\biggr]\\
&=\frac1t\int_0^{\infty} f_n'(v) e^{-tv}\td v\\
&=\dotsm\\
&=\frac1{t^{2n-1}}\int_0^{\infty} f_n^{(2n-1)}(v) e^{-tv}\td v\\
&=\frac{(-1)^n}{t^{2n-1}}\int_0^{\infty}\Biggl[\frac{1}{v}-\frac{1}{2}\coth\frac{v}{2} +\sum_{k=1}^n\frac{B_{2k}}{(2k)!} v^{2k-1}\Biggr]^{(2n-1)} e^{-tv}\td v\\
&=\frac{(-1)^n}{t^{2n-1}}\int_0^{\infty} K_{2n-1}(v) e^{-tv}\td v
\end{align*}
for $n\ge1$, where
\begin{equation}\label{K(n)-Bern-Coth-Eq}
K_{2n-1}(v)=\frac{\td^{2n-1}}{\td v^{2n-1}}\biggl(\frac{1}{v}-\frac{1}{2}\coth\frac{v}{2}\biggr)+\frac{B_{2n}}{2n}, \quad v>0, \quad n\ge1.
\end{equation}
\par
Using the relation
\begin{equation*}
\int_0^{\infty}\frac{x^{2n+1}\cos(ax)}{e^x-1}\td x=(-1)^n\frac{\td^{2n+1}}{\td a^{2n+1}}\biggl[\frac{\pi}{2}\coth(a\pi) -\frac{1}{2a}\biggr], \quad a>0, \quad n\ge0
\end{equation*}
in~\cite[p.~48, (D20)]{Flax-NASA-1970}, \cite[p.~506, 3.951.13]{Gradshteyn-Ryzhik-Table-8th}, and~\cite[p.~2]{Wintucky}, we arrive at
\begin{equation}\label{nth-deriv-mansour-K}
\frac{\td^{2n-1}}{\td v^{2n-1}}\biggl(\frac{1}{v}-\frac{1}{2}\coth\frac{v}{2}\biggr)
=(-1)^n2\int_0^{\infty}\frac{w^{2n-1}\cos(wv)}{e^{2\pi w}-1}\td w
\end{equation}
for $v>0$ and $n\ge1$.
Using the relation
\begin{equation*}
\int_0^{\infty}\frac{t^{2k-1}} {e^{2\pi t}-1}\td t=(-1)^{k-1}\frac{B_{2k}}{4k}, \quad k\geq 1
\end{equation*}
in~\cite[p.~29, (1.6.4)]{andrews-askey-roy-1999b}, \cite[p.~220]{Berndt}, and~\cite[p.~19]{Titchmarsh}, we obtain
\begin{equation}\label{nth-deriv-mansour-Bern}
(-1)^n2\int_0^{\infty}\frac{w^{2n-1}} {e^{2\pi w}-1}\td w=-\frac{B_{2n}}{2n}, \quad n\ge1.
\end{equation}
Substituting~\eqref{nth-deriv-mansour-K} and~\eqref{nth-deriv-mansour-Bern} into~\eqref{K(n)-Bern-Coth-Eq} reveals
\begin{equation*}
K_{2n-1}(v)=(-1)^{n-1}2\int_0^{\infty}\frac{w^{2n-1}[1-\cos(wv)]}{e^{2\pi w}-1}\td w , \quad v>0, \quad n\ge1.
\end{equation*}
Consequently, we have
\begin{equation}\label{R(n)-der-cos-int}
t^{2n-1}[-R_n'(t)]=2\int_0^{\infty} \biggl(\int_0^{\infty}\frac{w^{2n-1}[1-\cos(wv)]}{e^{2\pi w}-1} \td w\biggr)e^{-tv}\td v, \quad n\ge1.
\end{equation}
Hence, the function $t^{2n-1}[-R_n'(t)]$ is completely monotonic on $(0,\infty)$. This means that
\begin{equation}\label{Qi-Conj-Gam-Deg-der>(2n+1)}
\cmdeg{t}[-R_n'(t)]\ge2n-1, \quad n\ge1.
\end{equation}
\par
If the function $t^{\alpha}[-R_{n}'(t)]$ were completely monotonic on $(0,\infty)$, then its first derivative is negative, hence
\begin{equation*}
\alpha<-\frac{t[-R_{n}'(t)]'}{-R_{n}'(t)}=-\frac{tR_{n}''(t)}{R_{n}'(t)}.
\end{equation*}
From
\begin{equation*}
\lim_{t\to0^+}\biggl[t^2\biggl(\psi'(t)-\frac{1}{t}-\frac{1}{2t^2}\biggr)\biggr]
=\lim_{t\to0^+}\biggl[t^2\biggl(\psi'(t+1)-\frac{1}{t}+\frac{1}{2t^2}\biggr)\biggr]
=\frac12
\end{equation*}
and
\begin{equation*}
\lim_{t\to0^+}\biggl[t \biggl(\psi(t)-\ln t+\frac{1}{2 t}\biggr)\biggr]
=\lim_{t\to0^+}\biggl[t \biggl(\psi(t+1)-\ln t-\frac{1}{2 t}\biggr)\biggr]
=-\frac12,
\end{equation*}
it follows that
\begin{align*}
\lim_{t\to0^+}\biggl[-\frac{tR_{n}''(t)}{R_{n}'(t)}\biggr]
&=-\lim_{t\to0^+}\frac{t\bigl[\psi'(t)-\frac1t-\frac{1}{2t^2}-\sum_{k=1}^{n}\frac{B_{2k}}{t^{2k+1}}\bigr]} {\bigl[\psi(t)-\ln t +\frac{1}{2t}+\sum_{k=1}^{n}\frac{B_{2k}}{2k}\frac1{t^{2k}}\bigr]}\\
&=-\lim_{t\to0^+}\frac{t^2\bigl[\psi'(t)-\frac1t-\frac{1}{2t^2}\bigr]-\sum_{k=1}^{n}\frac{B_{2k}}{t^{2k-1}}} {t\bigl[\psi(t)-\ln t +\frac{1}{2t}\bigr]+\sum_{k=1}^{n}\frac1{2k}\frac{B_{2k}}{t^{2k-1}}}\\
&=2n.
\end{align*}
This means that
\begin{equation}\label{Qi-Mansour-degr<=2n}
\cmdeg{t}[-R_{n}'(t)]\le2n, \quad n\ge1.
\end{equation}
\par
Combining~\eqref{Qi-Mansour-degr<=2n} with~\eqref{Qi-Conj-Gam-Deg-der>(2n+1)} gives the inequality~\eqref{Qi-Conj-Gam-Deg-der(2n-1)-2n}.
The proof of Theorem~\ref{1der-deg=2n-1thm} is complete.

\section{Remarks}

Finally we list several remarks on main results and Qi's conjectures in this paper.

\begin{rem}
As in the derivation of~\eqref{K(n)-Bern-Coth-Eq}, integrating by parts $R_n'(t)$ inductively $2n\ge2$ times in~\eqref{Remainder-f-n-Rel-Eq} yields
\begin{equation*}
t^{2n}[-R_n'(t)]=\int_0^{\infty} K_{2n}(v) e^{-tv}\td v
\end{equation*}
for $n\ge1$, where, by~\eqref{id-gen-new-theta=1},
\begin{align*}
K_{2n}(v)&=\frac{\td^{2n}}{\td v^{2n}}\biggl(\frac{1}{v}-\frac{1}{2}\coth\frac{v}{2}\biggr)\\
&=\frac{\td^{2n}}{\td v^{2n}}\biggl(\frac{1}{v}-\frac{1}{2}-\frac{1}{e^v-1}\biggr)\\
&=\frac{(2n)!}{v^{2n+1}}-\sum_{p=1}^{2n+1}{(p-1)!S(2n+1,p)}\biggl(\frac1{e^v-1}\biggr)^p
\end{align*}
for $v>0$ and $n\ge1$.
It is easy to see that
\begin{align*}
K_2(v)&=\frac{2}{v^3}-\frac{2}{(e^v-1)^3}-\frac{3}{(e^v-1)^2}-\frac{1}{e^v-1}\\
&=\frac{2e^{3v}-e^{2v}(v^3+6)-e^v(v^3-6)-2}{(e^v-1)^3v^3},\\
\bigl[(e^v-1)^3v^3K_2(v)\bigr]'&=e^v\bigl[6e^{2v}-v^3-3v^2+6-e^v(2v^3+3v^2+12)\bigr]\\
&\to0,\quad v\to0^+,\\
\bigl(\bigl[(e^v-1)^3v^3K_2(v)\bigr]'e^{-v}\bigr)'&=12e^{2v}-e^v(2v^3+9v^2+6v+12)-3v(v+2)\\
&\to0,\quad v\to0^+,\\
\bigl(\bigl[(e^v-1)^3v^3K_2(v)\bigr]'e^{-v}\bigr)''&=24e^{2v}-e^v(2v^3+15v^2+24v+18)-6(v+1)\\
&\to0,\quad v\to0^+\\
\bigl(\bigl[(e^v-1)^3v^3K_2(v)\bigr]'e^{-v}\bigr)^{(3)}&=48e^{2v}-e^v(2v^3+21v^2+54v+42)-6\\
&\to0,\quad v\to0^+\\
\bigl(\bigl[(e^v-1)^3v^3K_2(v)\bigr]'e^{-v}\bigr)^{(4)}&=96e^v\biggl(e^v-1-v-\frac{9}{16}\frac{v^2}{2!}-\frac{1}{8}\frac{v^3}{3!}\biggr)\\
&>0
\end{align*}
for $v\in(0,\infty)$. Therefore, the function $K_2(v)$ is positive on $(0,\infty)$.
By Theorem~12b in~\cite[p.~161]{widder} stated in~\eqref{Theorem12b-Laplace}, we conclude that the function $t^2[-R_1'(t)]$ is completely monotonic on $(0,\infty)$, that is,
\begin{equation*}
\cmdeg{t}[-R_1'(t)]\ge2.
\end{equation*}
This supplies an alternative and partial proof for the second equality in~\eqref{Qi-Conj-Gam-Deg-der1-1-2}.
\end{rem}

\begin{rem}
Using~\eqref{R(n)-der-cos-int} and integrating by parts show
\begin{equation*}
t^{2n}[-R_n'(t)]=2\int_0^{\infty}\biggl[\int_0^{\infty}\frac{w^{2n}\sin(wv)}{e^{2\pi w}-1} \td w\biggr]e^{-tv}\td v.
\end{equation*}
When $n=1$, the second equality in~\eqref{Qi-Conj-Gam-Deg-der1-1-2} means that the function $t^2[-R_n'(t)]$ is completely monotonic on $(0,\infty)$. By Theorem~12b in~\cite[p.~161]{widder} stated in~\eqref{Theorem12b-Laplace}, we conclude
\begin{equation*}
\int_0^{\infty}\frac{u^2\sin(su)}{e^{u}-1} \td u\ge0, \quad s\in(0,\infty).
\end{equation*}
The equality~\eqref{Qi-Conj-Gam-Deg-der1(2n-1)} verified in Theorem~\eqref{1der-deg=2n-1thm} means that, when $n\ge2$, the function $t^{2n}[-R_n'(t)]$ is not completely monotonic on $(0,\infty)$. Again by Theorem~12b in~\cite[p.~161]{widder} stated in~\eqref{Theorem12b-Laplace}, we conclude that the function
\begin{equation}\label{u-n-power-sin<0}
\int_0^{\infty}\frac{u^{2n}\sin(su)}{e^{u}-1} \td u, \quad s\in(0,\infty) \quad n\ge2
\end{equation}
attains negative somewhere on $(0,\infty)$. Consequently, the function
\begin{equation*}
\int_0^{\infty}\frac{u[1-\cos(su)]}{e^u-1} \td u,\quad s\in(0,\infty)
\end{equation*}
is increasing on $(0,\infty)$ and the function
\begin{equation*}
\int_0^{\infty}\frac{u^{2n-1}[1-\cos(su)]}{e^u-1} \td u, \quad s\in(0,\infty), \quad n\ge2,
\end{equation*}
is decreasing on some subinterval $I\subset(0,\infty)$.
\par
From the formula~\eqref{u-n-power-sin<0} and
\begin{equation*}
\int_0^{\infty}\frac{x^{2n}\sin(ax)}{e^x-1}\td x=(-1)^n\frac{\td^{2n}}{\td a^{2n}}\biggl[\frac{\pi}{2}\coth(a\pi) -\frac{1}{2a}\biggr], \quad a>0, \quad n\ge0
\end{equation*}
in~\cite[p.~48, (D19)]{Flax-NASA-1970}, \cite[p.~506, 3.951.12]{Gradshteyn-Ryzhik-Table-8th}, and~\cite[p.~1]{Wintucky}, we reveal that the functions
\begin{equation*}
(-1)^n\biggl(\coth s-\frac{1}{s}\biggr)^{(2n)}\quad\text{and}\quad
(-1)^n\biggl(\frac{1}{e^s-1}-\frac{1}{s}\biggr)^{(2n)}
\end{equation*}
for $n\ge2$ attain negative on some subinterval $I\subset(0,\infty)$.
\end{rem}

\begin{rem}
From~\eqref{nth-deriv-mansour-K}, we conclude that
\begin{align*}
\int_0^{\infty}\frac{w^{2n-1}[1-\cos(wv)]}{e^{2\pi w}-1} \td w
&=\int_0^{\infty}\frac{w^{2n-1}}{e^{2\pi w}-1} \td w
-\int_0^{\infty}\frac{w^{2n-1}\cos(wv)}{e^{2\pi w}-1} \td w\\
&=\frac{1}{(2\pi)^{2n}}\int_0^{\infty}\frac{w^{2n-1}}{e^w-1} \td w
-\frac{(-1)^n}{2}\frac{\td^{2n-1}}{\td v^{2n-1}}\biggl(\frac{1}{v}-\frac{1}{2}\coth\frac{v}{2}\biggr)\\
&=\frac{\Gamma(2n)\zeta(2n)}{(2\pi)^{2n}}
-\frac{(-1)^n}{2}\frac{\td^{2n-1}}{\td v^{2n-1}}\biggl(\frac{1}{v}-\frac{1}{2}-\frac{1}{e^v-1}\biggr)\\
&=\frac{\Gamma(2n)\zeta(2n)}{(2\pi)^{2n}}
-\frac{(-1)^n}{2}\Biggl[\frac{(2n-1)!}{v^{2n}}-\sum_{p=1}^{2n}{(p-1)!S(2n,p)}\biggl(\frac1{e^v-1}\biggr)^p\Biggr],
\end{align*}
where $\zeta(z)$ denotes the Riemann zeta function~\cite[Chapter~25]{NIST-HB-2010}.
Consequently, we obtain three inequalities
\begin{align*}
\frac{(-1)^n}{2}\biggl(\frac{1}{v}-\frac{1}{2}\coth\frac{v}{2}\biggr)^{(2n-1)}
&<\frac{\Gamma(2n)\zeta(2n)}{(2\pi)^{2n}},\\
\frac{(-1)^n}{2}\biggl(\frac{1}{v}-\frac{1}{e^v-1}\biggr)^{(2n-1)}
&<\frac{\Gamma(2n)\zeta(2n)}{(2\pi)^{2n}},
\end{align*}
and
\begin{equation*}
\frac{(-1)^n}{2}\Biggl[\frac{(2n-1)!}{v^{2n}}-\sum_{p=1}^{2n}{(p-1)!S(2n,p)}\biggl(\frac1{e^v-1}\biggr)^p\Biggr]
<\frac{\Gamma(2n)\zeta(2n)}{(2\pi)^{2n}}
\end{equation*}
for $n\ge1$ and $v\in(0,\infty)$.
\par
Directly proving these three inequalities should be difficult and these functions involving in these three inequalities are important in special functions and number theory. These are two reasons why we write this remark.
\end{rem}

\begin{rem}
From the proof of Theorem~\ref{1der-deg=2n-1thm}, we can deduce that
\begin{align*}
\lim_{t\to\infty}R_n'(t)&=0, & \lim_{t\to\infty}\bigl[t^{2n-1}R_n'(t)\bigr]&=0,\\
\lim_{t\to\infty}\bigl[t^{2n+3}R_n'(t)\bigr]&=(-1)^{n+1}\frac{B_{2n+2}}{2n+2}, &
\lim_{t\to0^+}\bigl[t^{2n+1}R_n'(t)\bigr]&=(-1)^n\frac{B_{2n}}{2n}
\end{align*}
for $n\in\mathbb{N}$.
\end{rem}

\begin{rem}
The conjectures in~\eqref{R-n-0th-deg}, \eqref{R-01-m-th-deg}, and~\eqref{R-n-m-th-deg} posed by Qi in~\cite{CAM-D-18-02975.tex} are still kept open. Theorem~\ref{1der-deg=2n-1thm} does not give a full answer to the conjecture in~\eqref{Qi-Conj-Gam-Deg-der1(2n-1)}, but it still demonstrates that these open conjectures posed by Qi should be true.
\end{rem}

\begin{rem}
This paper is a revised version of electronic preprints~\cite{Qi-Mansour-Slovaca-5500.tex, Mansour-Qi-CMD.tex-arXiv, Mansour-Qi-CMD.tex-HAL, Qi-Mansour-Deg-Partial.tex}.
\end{rem}

\section{Acknowledgements and declarations}

\subsection{Acknowledgements}
The authors appreciate
\begin{enumerate}
\item
Ling-Xiong Han (Inner Mongolia University for Nationalities, China; hlx2980@163.com) for her careful computations and helpful discussions with the first author online via the Tencent QQ;
\item
Li Yin (Binzhou University, China; yinli7979@163.com) for his frequent communications and helpful discussions with the first author online via the Tencent QQ;
\item
anonymous referees for their careful corrections to and valuable comments on the original version of this paper.
\end{enumerate}

\subsection{Authors' contributions}
All authors contributed equally to the manuscript and read and approved the final manuscript.

\subsection{Conflict of interest}
The authors declare that they have no conflict of interest.

\subsection{Data availability statements}
Data sharing is not applicable to this article as no new data were created or analyzed in this study.

\end{document}